\documentclass[reqno, 12pt]{amsart}
\usepackage{amsfonts}
\usepackage{amssymb,amsmath}
\usepackage{amsthm}
\usepackage{upref}
\makeatletter
\def\LaTeX{\leavevmode L\raise.42ex
    \hbox{\kern-.3em\size{\sf@size}{0pt}\selectfont A}\kern-.15em\TeX}
 \makeatother

\numberwithin{equation}{section}

\newtheorem{lemma}{Lemma}[section]
\newtheorem{theorem}[lemma]{Theorem} 
\newtheorem{corollary}[lemma]{Corollary}
\newtheorem{proposition}[lemma]{Proposition}

\theoremstyle{definition}

\DeclareMathOperator{\Ran}{Ran}

\DeclareMathOperator{\spec}{spec}

\DeclareMathOperator*{\slim}{s-lim}

  \newcommand{\e}{\eqref}

\newcommand{\q}{\quad}

\newcommand{\ti}{\tilde}
\newcommand{\wt}{\widetilde}

\newcommand{\la}{\langle}
\newcommand{\ra}{\rangle}
\renewcommand{\d}{\delta}

   \newcommand{\sgn}{\operatorname{sgn}}

\renewcommand\Im{\operatorname{Im}}
\renewcommand\Re{\operatorname{Re}}

\newenvironment{pf}{\begin{proof}}{\end{proof}}

\def\qqq{\mathrel{\subset\mkern-15mu\lower.38ex\hbox{${\scriptscriptstyle\rightarrow}$}}}

\let\goth\mathfrak
\let\cal\mathcal

\let\Bbb\mathbb

   \DeclareMathOperator{\ac}{\rm ac}
  
\begin{document}

 \bigskip

\title[Carleman operator]
{Spectral and scattering theory for perturbations of the Carleman operator}
\author{ D. R. Yafaev}

 \dedicatory{To the memory of Vladimir Savel'evich Buslaev} 

 \address{ IRMAR, Universit\'{e} de Rennes I\\ Campus de
  Beaulieu, 35042 Rennes Cedex, FRANCE}
  
\email{yafaev@univ-rennes1.fr}



\begin{abstract}
We study spectral properties of the Carleman operator (the Hankel operator with kernel $h_{0}(t)=t^{-1}$)  and, in particular, find an explicit formula for its resolvent. Then we consider perturbations of the Carleman operator $H_{0}$ by Hankel operators $V$ with kernels $v(t)$
decaying sufficiently rapidly as $t\to\infty$ and not too singular at $t=0$. Our goal is to develop scattering theory for the pair $H_{0}$, $H=H_{0} +V $ and to construct an expansion in eigenfunctions of the continuous spectrum of the Hankel  operator $H$.   We also prove that under general assumptions  the singular continuous 
spectrum of the  operator $H$ is empty and that its eigenvalues may accumulate only to  the edge points $0$ and $\pi$ in  the
spectrum of $H_{0}$. We find  simple conditions for the finiteness   of the total number of eigenvalues of the operator $H$ lying above   the (continuous) spectrum of the Carleman operator $H_{0}$ and obtain an explicit estimate of this number. The theory constructed is somewhat analogous to the theory of one-dimensional differential operators.
\end{abstract}

\maketitle


\section{Introduction}  

{\bf 1.1.}
Hankel operators $H$ can be defined as integral operators,
\begin{equation}
(H f)(t) = \int_{0}^\infty h(t+s) f(s)ds, 
\label{eq:Hankel}\end{equation}
in the space $L^2 ({\Bbb R}_{+}) $ with kernels $h$ which depend  on the sum of variables only. As was pointed out by J.~S.~Howland in \cite{Howland}, self-adjoint Hankel  operators are to a certain extent similar to differential operators. In particular, Hankel  operators with continuous spectrum resemble singular differential operators.
In terms of this analogy, the Carleman operator $H_{0}$ corresponding to the kernel   $h_{0}(t)=t^{-1}$ plays the role of the ``free" Schr\"odinger operator $D^2$, $D=-i d /dx $,  in the space $L^2 ({\Bbb R}) $. The Carleman operator    can   easily be diagonalized by the Mellin transform.

   As far as the theory of Hankel operators is concerned, we refer to the books \cite{Pe} by V.~V.~Peller and \cite{Po}
by S.~R.~Power. We also note the paper
\cite{Howl0} by J.~S.~Howland where, in the trace class framework,  the structure of the absolutely continuous spectra of Hankel operators was described in terms of their symbols.

Our goal here is to study spectral  properties  of the Carleman operator $H_{0}=:{\bf C}$ and of Hankel operators $H$ with kernels $h(t)$ behaving asymptotically as $t^{-1}$ for $t\to \infty$ and $t\to 0$. In particular, we
develop the   scattering theory for  the pairs    $ H_{0} $, $H$. 
 
\medskip

{\bf 1.2.}
The first part of the paper (Sections~2, 3 and 4) is devoted to the study   of the Carleman operator $\mathbf{C}$ defined in the space $L^2 ({\Bbb R}_{+})$  by the relation 
\begin{equation}
(\mathbf{C} f)(t)=\int_{0}^\infty (t+s)^{-1} f (s) ds.
\label{eq:C}\end{equation}
It is easy to see that for all $k\geq  0$
\begin{equation}
 \int_{0}^\infty (t+s)^{-1}s^{-1/2\pm ik} ds= \lambda(k)t^{-1/2\pm ik} 
\label{eq:Ct}\end{equation} 
where
\begin{equation}
\lambda=\lambda(k)= \frac {\pi}{ \cosh (\pi k) } .
\label{eq:C6}\end{equation}
This equation establishes one-to-one correspondence between the quasi-momentum $k\geq 0$ and the energy $\lambda\in (0,\pi]$. It can be solved by the formula
\begin{equation}
k=k(\lambda) =\pi^{-1} \ln \big( (\pi + \sqrt{\pi^2-\lambda^2} )\lambda^{-1} \big) \geq 0.
\label{eq:B6k}\end{equation} 
Relations \e{eq:Ct} and \e{eq:C6} show that the spectrum of the operator $\bf C$ is absolutely continuous, coincides with the interval $[0,\pi]$ and has multiplicity two.

The dispersive relation \e{eq:C6}     plays    the same role for the Carleman operator as the relation $\lambda=k^2$ between the energy $\lambda>0$ and the  momentum $k>0$ for the differential  operator $D^2$ in the space $L^2({\Bbb R})$. In terms of this analogy, the eigenfunctions of the continuous spectrum $t^{-1/2\pm ik}$ of the operator $\bf C$ play the role of the eigenfunctions $e^{\pm ikx}$ of the operator $D^2$.  The singular points $t=0$ and $t=\infty$ correspond to the singular points $x=-\infty$ and $x=\infty$.

  One of our main results (Theorem~\ref{res}) yields the explicit formula for the resolvent ${\bf R}(z)= ({\bf C}-z I)^{-1}$ of the   operator $\mathbf{C}$. It plays the crucial role in our study of perturbations of the Carleman operator.
  
   In   Section~3 we study boundary values of ${\bf R}(z)$ as $z$ approaches the spectrum $[0,\pi]$ of the   operator $\mathbf{C}$ and, in particular, its edge points $z=\pi$ and $z= 0$. It turns out that the singularity of ${\bf R}(z)$  as $z\to\pi$ is quite similar to that of the resolvent $(D^2-zI)^{-1}$ as $z\to 0$. In particular, the   operator $\mathbf{C}$ has a resonance at the point $z=\pi$. On the contrary, the singularity of ${\bf R}(z)$  as $z\to 0$ is rather unusual and contains an oscillating factor.

  We also find in   Section~4   asymptotics of  the unitary group $\exp(-i {\bf C}T) $   as $T\to \pm \infty$. We show that functions $(\exp(-i {\bf C}T)f)(t)$ are localized for large $|T|$ in exponentially small neighbourhoods of singular points $t=0$ and $t=\infty$. This should be compared with the well known fact that functions  $(\exp(-i D^2 T)f)(x)$ ``live" in the region where $|x|$ and $|T|$ are of the same order. Our formulas for $\exp(-i {\bf C}T) $ (see   Theorem~\ref{Time})  are somewhat similar but
   essentially more complicated than those for the unitary group $\exp(- i D^2 T   ) $.

 \medskip

{\bf 1.3.}
In the second part of the paper (Sections~5 and 6) we study perturbations  of the Carleman operator $H_{0}= {\bf C}$  by Hankel  operators $V$,
\begin{equation}
(V f(t)=\int_{0}^\infty v(t+s)  f (s) ds,
\label{eq:H1}\end{equation}
 with kernels $v(t)$ decaying faster that $t^{-1}$ as $t\to \infty$ and less singular than $t^{-1}$ as $t\to 0$. Such perturbations of ${\bf C}$
play the role of perturbations of the operator $D^2$ by operators of the multiplication by functions ${\sf V}(x)$ decaying sufficiently rapidly as $|x|\to \infty$.
 To a certain extent, Section~5 can be considered as a continuation of the paper \cite{Howland} 
 where the Mourre method was used for the study of Hankel operators. Nevertheless, our approach to this problem is quite different from that of  \cite{Howland}. Actually, we   follow the analogy with the one-dimensional Schr\"odinger operator (see, e.g., the original paper \cite{F} by L.~D.~Faddeev or the book \cite{Yafaev1}). 
 
  There is, however, an important difference between 
   one-dimensional Schr\"odinger  operators  and Hankel operators. In the first case one can develop the theory relying exclusively on Volterra integral equations while such a possibility is of course lacking in the second case. In this respect, the theory
   of Hankel operators is closer  to  the theory
   of    one-dimensional differential operators of   order higher than two where Fredholm integral equations occur naturally  (see \cite{LNM}). Note that for the operator $D^n$ in $L^2 ({\Bbb R})$ the eigenfunctions $e^{\pm ikx}$ are   the same as those  for the operator $D^2$ for all $n=1,2,\ldots$, but the dispersive relation has the form $\lambda=k^n$. For the Carleman operator, the dispersive relation \e{eq:C6} is more complicated.
   
In Section~5,   we proceed from the results of Section~3 on the existence of suitable boundary values of the resolvent $R_{0}(z)=(H_{0}-zI)^{-1}$ (the limiting absorption principle for the operator $H_{0}$). Then we apply general results of abstract scattering theory (see, e.g., the paper \cite{Kuroda} or the book \cite{Yafaev}) and establish the limiting absorption principle for the operator $H=H_{0}+V$.
    This allows us   to obtain rather a detailed information about eigenfunctions of the continuous spectrum of the operator $H $ and to prove
  an expansion in eigenfunctions of  $H$. We then find expressions   for the wave operators for the pair $H_{0}$, $H$ and for  the corresponding scattering matrix in terms of the asymptotics of the eigenfunctions of $H$  as $t\to \infty$ and $t\to 0$. 
  Formulas obtained look similar to those for   one-dimensional differential operators.

 Finally, in Section~6 we study the discrete spectrum of the operator $H=H_{0}+V$ lying above its continuous spectrum $[0,\pi]$. We show that it consists only of a finite number of eigenvalues  if the function $v (t)$ decays sufficiently rapidly as $t\to \infty$ and it is not too singular as $t\to 0$. Moreover, the operator $H $ necessary has an eigenvalue larger than $\pi$ if $V\geq 0$ and $V\neq 0$. These results are    similar to those on the negative spectrum of the Schr\"odinger operator $D^2 +{\sf V}(x)$. This is of course quite   natural because the singularities of the resolvents at the corresponding edge points of the continuous spectra are also similar.

 On the contrary, the finiteness of the negative spectrum of $H=H_{0}+V$ is not determined by the behaviour of $v (t)$ at singular points $t=0$ and $t=\infty$. For example, one can construct functions $\eta(x)$ from the Schwartz class (in fact,   the Fourier transforms of $\eta$ even belong to the class $C_{0}^\infty({\Bbb R})$) such that   for the kernel $v(t)=t^{-1} \eta(\ln t)  $, the negative spectrum
   of $H$ is infinite. This phenomenon is of course related to a complicated structure of the ``free" resolvent $R_{0} (z)$ as $z\to 0$.  But this is a subject of another paper.
   
Actually, in Sections~5 and 6 we admit perturbations   by sufficiently general  integral operators $V$ (not necessarily Hankel  operators) with kernels $ {\bf v}(t,s)$ satisfying some  decay assumptions    at infinity and some regularity assumptions at  the origin.

\medskip

{\bf 1.4.}
As was already mentioned, spectral and scattering theories of differential operators and of Hankel operators are to a certain extent parallel.  It might be interesting to extend this analogy a bit further. As examples, let us mention various trace formulas (cf. the paper \cite{BF} by V. S. Buslaev and L. D. Faddeev where the Schr\"odinger operator  on the half-axis was considered) and the inverse problem of a reconstruction of the kernel given the corresponding scattering matrix  (cf. the paper \cite{F} by  L. D. Faddeev where the Schr\"odinger operator  on the whole axis was considered). Note, however, that such more advanced questions might be difficult already for differential operators of order higher than two and even more difficult for Hankel operators. In this respect, we mention the paper \cite{Ost} where a trace formula (to put it differently, an expression for the perturbation determinant in terms of solutions of the corresponding differential equation) was obtained for  differential operators of an arbitrary order. The solution of the inverse problem (see, e.g., the book \cite{BDT})  for  differential operators of an arbitrary order seems to be in a   less satisfactory state than for order  two.

 \medskip
 
 The author is grateful for hospitality and financial
support to the  Mittag-Leffler Institute, Sweden, where the paper was completed 
during the author's stay in the autumn of 2012.

\section{The resolvent of the Carleman operator}  

Here we calculate the resolvent ${\bf R}(z)= ({\bf C}-zI)^{-1}$ of the Carleman operator $\mathbf{C}$ defined in the space $L^2 ({\Bbb R}_{+})$  by   relation 
\e{eq:C}.

\medskip

{\bf 2.1.}
Let us introduce the Mellin transform 
$M$,
\begin{equation}
  ( M f)(k)=  (2\pi)^{-1/2  }\int_{0}^\infty t^{-1/2 -ik } f (t) dt,
\label{eq:C3}\end{equation}
which is the    unitary mapping $M :L^2 ({\Bbb R}_{+})\to L^2 ({\Bbb R} )$.
Obviously, we have 
\begin{equation}
(  M  {\bf C} f )(k)=  (2\pi)^{-1/2  }\int_{0}^\infty ds f(s)\int_{0}^\infty dt t^{-1/2-ik} (t+s)^{-1}  = \lambda(k) ( M f )(k)
\label{eq:C4}\end{equation}
where
\begin{equation}
\lambda(k)=  \int_{0}^\infty t^{-1/2 - ik}  (t +1)^{-1} dt .
\label{eq:C5}\end{equation}
As is well known and as we shall see below, this integral is given by formula \e{eq:C6}.
Thus the spectrum of the  operator ${\bf C}$ is absolutely continuous, it coincides with the interval $[0,\pi]$ and has multiplicity $2$. 
 
To calculate the resolvent  of the   operator ${\bf C}$, we have to solve the equation
\begin{equation}
({\bf C} -zI) f = f_{0},\q z\in{\Bbb C}\setminus [0, \pi].
\label{eq:C7}\end{equation}
Making here the Mellin transform
and using  relation \e{eq:C4}, we can equivalently rewrite \e{eq:C7} as the equation
\[
(\lambda (k) -z) \ti{f}(k)=    \ti{f}_{0}(k) 
\]
for the functions $  \ti{f}_{0}= M f_{0}$ and  $   \ti{f} = M  f $. In view of   \e{eq:C6},    this equation can be solved by the formula
\begin{equation}
 \ti{f}(k)=(\lambda (k) -z)^{-1}  \ti{f}_{0}(k) =   -z^{-1}  \ti{f}_{0}(k)  -\pi z^{-2} \frac {1}{ \cosh (\pi k)-\pi z^{-1} }
   \ti{f}_{0}(k) .
\label{eq:C9}\end{equation}

Now we have to make here the inverse Mellin transform. To that end, we use the following elementary assertion.

\begin{lemma}\label{mellin}
Let $Q$ be the operator of multiplication in the space $L^2 ({\Bbb R})$ by a  function $q\in L^\infty (\Bbb R)\cap L^1 (\Bbb R)$. Then 
\begin{equation}
(M^* Q M f) (t)=\int_{0}^\infty (ts)^{-1/2} {\bf q}(t/s) f(s)ds
\label{eq:C8A}\end{equation}
  where
\begin{equation}
{\bf q}(u)= (2\pi)^{-1} \int_{-\infty}^\infty u^{ik} q(k)dk.
\label{eq:C8B}\end{equation}
 \end{lemma}
 
 \begin{pf}
 Interchanging the order of integrations, we see that
 \begin{align*}
 (M^* Q M f)(t)=&(2\pi)^{-1} \int_{- \infty}^\infty dk t^{-1/2+ik} q (k)
 \int_{0}^\infty ds f(s) s^{-1/2-ik}
 \\
 =&(2\pi)^{-1} \int_{0}^\infty ds f(s)
  \Big(\int_{- \infty}^\infty  t^{-1/2+ik} s^{-1/2-ik} q(k)dk\Big).
 \end{align*}
 This yields formulas \e{eq:C8A}, \e{eq:C8B}.
   \end{pf}
   
 This result shows that if  $\phi$ is a bounded  function of $\lambda\in[0,\pi]$ such that $\lambda^{-1}\phi(\lambda)$ belongs to $L^1 (0,\pi)$, then 
 $\phi (\bf C)$  is an integral operator acting by formula \e{eq:C8A}. Here  the function ${\bf q}(u)$ is defined by equality \e{eq:C8B}  with $q(k)= \phi (\frac{\pi}{\cosh(\pi k)})$.
 
 In view of Lemma~\ref{mellin} it follows from relation \e{eq:C9} that
\begin{equation}
({\bf R}(z)  f_{0})(t)= f(t)=   -z^{-1}  f_{0}(t)  -  z^{-2}\pi^{-1} \int_{0}^\infty (ts)^{-1/2}  {\cal I}_{z}(t/s) f_{0}(s) ds
\label{eq:C11}\end{equation}
where
\begin{equation}
 {\cal I}_{z}(u)=   \frac{\pi}{2} \int_{- \infty}^\infty    \frac {u^{ik}}{ \cosh (\pi k)-\pi z^{-1} }  dk.
\label{eq:C12}\end{equation}

\medskip

{\bf 2.2.}
Let us calculate integral \e{eq:C12}. We set $\zeta= \pi z^{-1}$ so that $\zeta\in{\Bbb C}\setminus [1,\infty)$ if $z \in{\Bbb C}\setminus [0, \pi]$. Making the change of variables $p =e^{\pi k}$ we see that  
\begin{equation}
 {\cal I}_{z} (u) =    \int_{0}^\infty    \frac {p^{i x}}{ p^2-2 \zeta p +1}  dp,\q x=\pi^{-1}\ln u.
\label{eq:B}\end{equation}

We are going to calculate this integral by residues.
  Observe that the equation $p^2-2 \zeta p +1=0$ has two roots
\begin{equation}
p_{1}(\zeta)=  \zeta + \sqrt{\zeta^2-1}, \q p_{2} (\zeta)= \zeta -  \sqrt{\zeta^2-1}, 
\label{eq:B1}\end{equation}
which are different if $\zeta \neq -1$ and $ p_1 (\zeta)  p_2 (\zeta)= 1$.  
We fix $\arg p$ in the complex plane with the cut along $[0,\infty)$ by the condition $\arg p\in [0, 2\pi]$. Then
\[
  \arg p_1 (\zeta)+  \arg p_2 (\zeta)= 2\pi,
\]
and hence
\begin{equation}
\ln p_1 (\zeta)+ \ln p_2 (\zeta) = 2\pi i.
\label{eq:Bzz}\end{equation}

   Let us consider the contour $C_{\varrho}$ in ${\Bbb C}\setminus [0,\infty)$ consisting of the interval $[0, \varrho]$ on the upper edge of the cut, the circle $|p |= \varrho$ and the interval $[\varrho, 0]$ on the lower edge of the cut. By the Cauchy theorem for  sufficiently large $\varrho$, we have
\begin{equation}
   \int_{C_{\varrho}}     \frac {p^{i x}}{ p^2-2 \zeta p +1}  dp
=2 \pi i  \sum_{j=1} ^2{\rm Res}_{p=p_{j}(\zeta)}   \frac {p^{i  x}}{ p^2-2 \zeta p +1} .
\label{eq:B2}\end{equation}
Computing the residues, we find that the right-hand side here equals
\begin{equation}
   2 \pi i   \frac {p_1(\zeta)^{i  x} - p_2(\zeta)^{i  x}}{ p_1(\zeta) - p_2(\zeta)}=\frac{\pi i}{\sqrt{\zeta^2 - 1}}\big( p_1(\zeta) ^{ix}- p_2(\zeta)^{ix }\big) .
\label{eq:B3}\end{equation}
Note that this expression does not depend on the choice of the sign of
 $\sqrt{\zeta^2 - 1}$.
In the left-hand side of \e{eq:B2},  the integral over the lower edge of the cut equals 
\[
\int_{\varrho}^0  \frac {(pe^{2\pi i})^{i x}}{ p^2-2\zeta p +1}  dp
=- e^{-2\pi x  } \int_0^{\varrho}   \frac {p^{i x}}{ p^2-2 \zeta p +1}  dp.
\]
The integral over the circle $|p|= \varrho$ tends to zero as $\varrho \to \infty$ because   $| p^{i x}|=e^{-x\arg p }$ is bounded by $1$ for $x\geq 0$ and by $e^{- 2\pi x }$ for $x <0$. Therefore passing in \e{eq:B2} to the limit $\varrho \to \infty$, we obtain the equation for integral \e{eq:B}:
\begin{equation}
(1-u^{-2   }) {\cal I}_{z} (u)=  \frac{\pi i}{\sqrt{\zeta^2 - 1}}\big(  p_1(\zeta) ^{i x}- p_2(\zeta)^{i x}\big), \q \zeta=\pi/z .
\label{eq:B4x}\end{equation}

If $\zeta=-1$, then $p_1(\zeta)=p_2(\zeta)= -1$ so that the right-hand sides of formulas \e{eq:B2} and hence of \e{eq:B4x} should be replaced by 
\[
2 \pi i   \, {\rm Res}_{p=-1}   \frac {p^{i x}}{ (p +1)^2}= 2 \pi i  \frac {d}{ dp} p^{i x}\big|_{p=-1}= -2 \pi x (e^{\pi i})^{ix- 1}=2 \pi x  e^{-\pi x} .
\]
It follows that
\[
  {\cal I}_{-\pi} (u) =  \frac{ 2  \ln u }{ u-u^{-1} }  .
\]

Let us formulate the result obtained.

\begin{lemma}\label{int}
Integral \e{eq:C12} for $z \in{\Bbb C}\setminus [0,\pi]$ is given by the equation 
\begin{equation}
(1-u^{-2  })  {\cal I}_{z} (u)=  \frac{\pi  i}{\sqrt{\zeta^2 - 1}}
\Big(  u ^{i /\pi \ln p_{1} (\zeta) } 
-   u ^{ i  /\pi \ln p_{2} (\zeta) } \Big), \q \zeta=\pi/z ,
\label{eq:B4}\end{equation}
where  the numbers $p_{j}(\zeta)$ are defined by formulas  \e{eq:B1} and $\arg p_{j}(\zeta)\in (0,2\pi)$.
\end{lemma}

 Consider the particular case $\zeta=0$. If we, for example,  choose $\sqrt{\zeta^2-1}=i$, then $\ln p_{1}(\zeta)=\pi i/2$, $\ln p_{2}(\zeta)=3\pi i/2$ and hence
   formula \e{eq:B4} for integral \e{eq:C12} yields
   \[
   \int_{0}^\infty \frac{p^{ix}} {p^2+1} dp= \frac{\pi}{2\cosh(\pi x/2)}.
   \]
   This is equivalent to    expression \e{eq:C6}  for integral \e{eq:C5}.

  \medskip  

{\bf 2.3.}
To state the formula  for  the resolvent of the Carleman operator, we first rewrite equation \e{eq:B4}  in terms of the variable $z=\pi\zeta^{-1}$. Let us consider the function
\begin{equation}
\varphi(z)= \sqrt{z^2 - \pi^2}
\label{eq:B4z}\end{equation}
 in the complex plane cut along $[-\pi,\pi]$ and fix its branch by the condition $\varphi(z)>0$ for $z>\pi$.  
Observe that the function 
\begin{equation}
q(z)=\frac{\pi-i  \sqrt{z^2 - \pi^2}}{z},\q z\in{\Bbb C}\setminus [-\pi,\pi],
\label{eq:B4zz}\end{equation}
does not take positive values which allows us to set $\arg q(z)\in (0,2\pi)$.  With this convention, the function 
\begin{equation}
{\bf k}(z)=\frac{1}{\pi} \ln q(z)
\label{eq:kk}\end{equation}
is   analytic   for $z\in{\Bbb C}\setminus [-\pi,\pi]$. Since
$\sqrt{\zeta^2 - 1}= i z^{-1}\sqrt{z^2 - \pi^2}$ (this fixes the sign of the left-hand side), we have $  p_{2}(\pi/z) = q(z)$ and $\ln p_{2}(\pi/z)=\pi {\bf k}(z)$. 
Using \e{eq:Bzz}, we  also see that
$\ln p_1 (\pi/z)=- \pi {\bf k} (z) +2\pi i$, and hence formula \e{eq:B4} can be rewritten as
 \begin{equation}
(1-u^{-2  }) {\cal I}_{z} (u)=  \frac{\pi z}{\sqrt{z^2 - \pi^2}}
\big( u^{-2} u ^{-i {\bf k} (z)  } 
-   u ^{ i {\bf k} (z)    } \big). 
\label{eq:B4aa}\end{equation}

  Putting together formulas \e{eq:C11} and \e{eq:B4aa}, we obtain the expression for  the resolvent of the Carleman operator.

\begin{theorem}\label{res}
Let the function ${\bf k} (z)   $ be defined for $z\in {\Bbb C}\setminus [-\pi,\pi]$  by    formulas \e{eq:B4zz},  \e{eq:kk} and the condition $\arg q(z)\in (0,2\pi)$. Set
\begin{equation}
 \rho  (u;z)=  \frac{  u ^{i {\bf k}(z) }} {u^{-2  }-1}
+ \frac{  u ^{-i {\bf k}(z) }} { u^{2  } -1}   .
\label{eq:B4xx}\end{equation}
Then the resolvent     ${\mathbf{R}}(z) =({\mathbf{C}}-z I)^{-1}$ of operator \e{eq:C}  admits the representation
\begin{equation}
 {\mathbf{R}} (z)=-z^{-1}( I +   {\mathbf{A}} (z)) 
\label{eq:B5}\end{equation}
where $ {\mathbf{A}} (z) $ is the integral operator with kernel  
\begin{equation}
{\bf a} (t,s; z)= \frac{ 1}{\sqrt{z^2 - \pi^2}} (ts)^{-1/2} \rho (t/s;z).
\label{eq:B5y}\end{equation}
 \end{theorem}
 
 Lemma~\ref{api} shows that Theorem~\ref{res} remains true for all $z\in {\Bbb C}\setminus [0,\pi]$.
 
 Let us now discuss properties of the function $ \rho  (u;z)$. We first observe that function \e{eq:B4z} satisfies the condition
$\varphi(\bar{z})=\overline{\varphi(z)}$,   $\varphi (z)<0$ for $z<-\pi$ and
 \begin{equation}
\varphi(\lambdaÊ\pm i0) = \pm i\sqrt{\pi^2-\lambda^2}, \q \lambda\in [-\pi, \pi].
\label{eq:xr1}\end{equation}
In the following assertion we collect necessary properties of the function $ {\bf k}(z)$.

\begin{lemma}\label{k}
$1^0$ The function ${\bf k} (z) $ is an analytic function of $z\in {\Bbb C}\setminus [-\pi,\pi]$ and it satisfies the identity
\begin{equation}
{\bf k}(\bar{z})=- \overline{{\bf k}(z)}.
\label{eq:Bz}\end{equation}

$2^0$ The limits of ${\bf k} (z)   $ on the cut exist $($except the point $z=0)$ and
\begin{equation}
{\bf k}(\lambda+i0)=  k(|\lambda|) + i ,\q \lambda\in [-\pi,0),
\label{eq:Kk}\end{equation}
\begin{equation}
{\bf k}(\lambda+i0)=  k(\lambda) + 2 i ,\q \lambda\in (0,\pi],
\label{eq:Kk1}\end{equation}
where $k(\lambda) $ is function \e{eq:B6k}.
 \end{lemma}

\begin{pf}
Let $q (z)   $ be function \e{eq:B4zz}.
 If $z>\pi$, then  \e{eq:Bz} is true because $|q(z)|=1$ and hence $\Re {\bf k}(z)=0$. Then, by analytic continuation, \e{eq:Bz} extends to all complex $z$.

  It follows from \e{eq:xr1} that
\begin{equation}
q(\lambda+i 0)= (\pi +  \sqrt{  \pi^2 -\lambda^2}) \lambda^{-1}.
\label{eq:Kk2}\end{equation}
Therefore  $  q(\lambda+i0)< 0$ and hence $\arg q(\lambda+i0)=\pi$ for $\lambda\in [-\pi,0)$, which  proves \e{eq:Kk}. If  $\lambda\in (0,\pi]$, then $  q(\lambda+i0)>0$. To calculate
$\arg q(\lambda+i0) $, we  pass from the half-line $\lambda>\pi$ to the upper edge of the cut $(0,\pi)$ around the point $\lambda=\pi$ by a small semi-circle lying in the upper half-plane. Observe that $\Im q(\lambda)<0$ for $\lambda>\pi$ and that $q(z)$ arrives to the positive value \e{eq:Kk2}  remaining always in the lower half-plane. Therefore $\arg q(\lambda+i0)=2\pi$, which proves \e{eq:Kk1}.    
\end{pf}

Let us come back to the function $ \rho  (u;z)$.

\begin{proposition}\label{rho}
Let the function $\rho  (u;z) $ be defined  for $u>0$ and    $z\in {\Bbb C}\setminus [-\pi,\pi]$  by    formula \e{eq:B4xx}. Then:

$1^0$
The function $\rho  (u;z) $ depends analytically on  $z\in {\Bbb C}\setminus [-\pi,\pi]$, and it  is a $C^\infty$ function of $u\in{\Bbb R}_{+}$; in particular,  we have 
 \[
\rho  (1; z)= -  1- i {\bf k} (z).
\]  

$2^0$
The function $\rho  (u;z) $ satisfies the identities 
 \begin{equation}
\rho  (u^{-1}; z)= \rho  (u;z)
 \label{eq:Kk3}\end{equation}
 and
 \begin{equation}
\rho  (u ; \bar{z}  )= \overline{\rho  (u; z)}.
 \label{eq:Bzg}\end{equation}

 $3^0$ 
  Set   $\kappa(z)=\pi^{-1}\min\{\arg q(z),2\pi - \arg q(z)\} >0$  for $z\not\in [0,\pi]$ $($note that $\kappa(z)=1$ for  $z \in [-\pi,0))$. Then  $\rho  (u;z)= O(u^{-\kappa(z) })$ as $u\to \infty$, 
$\rho  (u;z)= O(u^{\kappa(z) })$ as $u\to 0$.
 
 $4^0$
 For all $u> 0$,
the limits of $\rho  (u;z)  $ on the cut exist $($except the point $z=0)$ and
\begin{align}
\rho (u; \lambda+i0)&=\frac{u^{i k(|\lambda|)} - u^{-i k(|\lambda|)}}{u^{-1} -u} ,\q \lambda\in [-\pi,0),
\label{eq:rr1}\\
 \rho  (u;\lambda+ i0 )&=  \frac{  u ^{i k(\lambda) }} {1-u^{2  } }
+ \frac{  u ^{-i k(\lambda) }} {1- u^{-2  }  }  ,\q \lambda\in (0,\pi].
\label{eq:B4xxx}\end{align}

 $5^0$ 
  The function $\rho (u; z)$ is uniformly bounded in $u\in{\Bbb R}_{+}$ and $z$ belonging  to compact subsets of $\Bbb C\setminus\{0\}$ including the values of $z$ on the cut along $[-\pi,\pi]$.
 \end{proposition}

\begin{pf}
Statement $1^0$ and identity \e{eq:Kk3} are obvious.
Identity \e{eq:Bzg} follows from \e{eq:Bz}.  
It follows from formulas \e{eq:kk} and \e{eq:B4xx}  that
\[
|\rho  (u; z)|\leq \big|\frac{  u ^{-\pi^{-1} \arg q(z)}  } {u^{-2  }-1} \big|
+  \big|\frac{  u ^{\pi^{-1} \arg q(z) }  } {u^{2  }-1}\big| .
\]
This immediately implies statement $3^0$. Representations \e{eq:rr1} and \e{eq:B4xxx} are direct consequences of equalities \e{eq:Kk} and \e{eq:Kk1}, respectively.  Statement $5^0$ is a direct consequence of definition \e{eq:B4xx} (see also formulas \e{eq:rr1} and \e{eq:B4xxx}).
\end{pf}

   Various properties of the resolvent ${\bf R} (z)$ are direct consequences of Proposition~\ref{rho}. Putting together identities \e{eq:Kk3} and \e{eq:Bzg}, we see that
   $\rho  (u ; \bar{z}  )= \overline{\rho  (u^{-1}; z)}$ and hence $\overline{a(t,s;z)}= a(s,t; \bar{z})$   which is consistent with the relation $\mathbf{R} (\bar{z})=\mathbf{R}^* (z)$ (the self-adjointness  of $\bf C$).
   
   Let us define the complex conjugation ${\cal C}$   by the equality 
     \begin{equation}
     ({\cal C}f)(t)=\overline{f(t)}.
     \label{eq:CC}\end{equation}  
     Identity \e{eq:Bzg} shows that $\overline{a(t,s;z)}= a(t,s; \bar{z})$ which is consistent with the relation ${\cal C} R( z)=R(\bar{z}){\cal C}$ (the invariance of $\mathbf{C}$ with respect to the complex conjugation). Thus we have $a(t,s;z)=a(s, t;z)$ which is actually a consequence of \e{eq:Kk3} solely.

 The function $\rho(u;z)$ is   analytic for $z\in{\Bbb C}\setminus [-\pi,\pi]$ only while the resolvent ${\bf R}(z)$ and hence ${\bf A}(z)$ should be analytic for all $z\in{\Bbb C}\setminus [0,\pi]$. To see this fact directly, we have to observe that the limits ${\bf a}(t,s;\lambda\pm  i0)$ exist and 
  \[
{\bf a}(t,s;\lambda +  i0) = {\bf a}(t,s;\lambda - i0), \q \lambda\in [-\pi, 0).
\]
Indeed, combining identities   \e{eq:Bzg} and \e{eq:rr1}, we see that
 $\rho (\lambda+i0)=- \rho (\lambda-i0)$.  So it remains to take equality   \e{eq:xr1} into account. Let us state the result obtained.

\begin{lemma}\label{api}
The function ${\bf a}(t,s; z)$ is  analytic for  $z\in{\Bbb C}\setminus [0,\pi]$
 and
 \begin{equation}
{\bf a} (t,s; \lambda)= \frac{ 2  }{\sqrt{ \pi^2-\lambda^2}}
\frac {\sqrt{ ts}}{s^2-t^2 }
\sin\big(k(|\lambda|) \ln( t/s)\big) ,\q \lambda\in (-\pi,0),
\label{eq:RR}\end{equation}
where   $k(|\lambda|)$ is determined by equality  \e{eq:B6k}. In particular, we have
\[
{\bf a} (t,s; -\pi)= \frac{ 2  }{  \pi^2 }
\frac {\sqrt{ ts}}{s^2-t^2 } \ln( t/s).
\]
 \end{lemma}

 Finally, we note that according to parts $1^0$ and $3^0$ of Proposition~\ref{rho}
\[
\int_{0}^\infty |\rho  (u;z)| u^{-1} du<\infty.
\]
  This estimate allows one to give a direct proof that the integral operator $ {\mathbf{A}} (z) $ with kernel \e{eq:B5y}  is bounded.

\section{Approaching the continuous spectrum}  

Now we discuss   boundary values of the resolvent ${\bf R} (z)= ({\bf C}-z I)^{-1}$ as $z$ approaches the cut along $[0,\pi]$.

\medskip

{\bf 3.1.}
Recall that  the kernel   of the operator $\mathbf{A}(z)$ is given by formula \e{eq:B5y}. 
It follows from part $4^0$ of Proposition~\ref{rho} (see also equality \e{eq:xr1}) that for all $t,s >0$ there exists the limits 
 \begin{equation}
{\bf a} (t,s; \lambda+i 0)= -  i\frac{1}{ \sqrt{\pi^2 - \lambda^2}} (ts)^{-1/2} \rho (t/s; \lambda+i0), \q \lambda\in (0,\pi),
\label{eq:B5yy}\end{equation}
where the function $\rho (t/s; \lambda+i0) $ is given by formula \e{eq:B4xxx}.

Let us consider    the spectral family ${\bf E}(\lambda)$  of the operator $\bf C$. Using that
\[
2\pi i ( d {\bf E} (\lambda) f, g)/ d\lambda =( {\bf R} (\lambda+i0)f,g) - ({\bf R}(\lambda-i0)f,g),
\]
it is easy to derive from formulas  \e{eq:B5}, \e{eq:B4xxx}  and \e{eq:B5yy} the following result.

\begin{lemma}\label{ee}
For all $t,s >0$,
the integral   kernel $e(t,s; \lambda)$ of the operator ${\bf E}(\lambda)$  is differentiable in $\lambda$ and
\[
\frac{e(t,s; \lambda)}{d\lambda}=\frac{1}{\pi\lambda \sqrt{\pi^2 - \lambda^2}} (ts)^{-1/2}\cos \big(k(\lambda) \ln (t/s)\big),\q \lambda\in (0,\pi).
\]
 \end{lemma} 
 
 \medskip

{\bf  3.2.}
Here we collect results that will be used in Section~5. Let   the operator $Q$ in the space $L^2 ({\Bbb R}_{+})$ be defined by the equality
\begin{equation}
(Q f)(t)= \la \ln t \ra f(t)\q {\rm where} \q \la \ln t \ra=(1+ |\ln t|^2)^{1/2}.
\label{eq:QQ}\end{equation}
According to formula \e{eq:B5y} it follows from
parts $4^0$ and $5^0$ of Proposition~\ref{rho}
  that the kernel of the  operator 
$Q^{-\beta} {\bf A}(z) Q^{-\beta}$ depends continuously on $z$,  and it is uniformly bounded
by\footnote{Here and in what follows we denote by $C$ (with various indices) positive constants whose values are of no importance.}
\[
C (ts)^{-1/ 2}  \la \ln t \ra^{-\beta}  \la \ln s \ra^{-\beta}.
\]
This function belongs to $L^2 ({\Bbb R}_{+}\times {\Bbb R}_{+})$ for $\beta>1/2$. Therefore by  the Lebesgue dominated convergence theorem, the
    operator-valued function
$Q^{-\beta} {\bf A}(z) Q^{-\beta}$   depends continuously on $z$ in the Hilbert-Schmidt norm. In view of  representation \e{eq:B5} this leads to the following assertion.

\begin{proposition}\label{cont}
For all $\beta>1/2$,  the operator-valued function
$Q^{-\beta} {\bf R}(z) Q^{-\beta}$ depends continuously in the norm of the space $L^2 ({\Bbb R}_{+} )$ on $z$ in the complex plane cut along $[0,\pi]$   as $z$ approaches the cut with exception of the points $0$ and $\pi$. Moreover, this function is H\"older continuous with any exponent $\gamma< \beta -1/2$
$($and $\gamma\leq 1)$.
 \end{proposition}
 
In formulas below $\lambda$ and $k$ are related   by equalities \e{eq:C6} or \e{eq:B6k}.     It follows from relation  \e{eq:B4xxx} that
for all $\lambda\in (0,\pi)$ there exist
\[
 \lim_{  u\to  \infty     }u^{ ik} \rho  (u; \lambda+ i0)= \lim_{  u \to 0   }u^{- ik} \rho  (u; \lambda+ i0)=  1. 
\]
 Therefore, again by the Lebesgue dominated convergence theorem, formula \e{eq:B5yy} yields the asymptotics of the function $({\bf A} (\lambda+ i0) f)(t)$ as $t\to \infty$ and as $t\to 0$. 
 Taking also into account representation \e{eq:B5}, we can state the following result.

    \begin{proposition}\label{asympt1x}
Suppose that 
\begin{equation}
\int_{0}^\infty t^{-1/2} |f(t)| dt <\infty
\label{eq:Li2x}\end{equation}
and   that $f(t)=  o (t^{-1/2})$ as $t \to \infty$ and as $t\to 0$.
 Then for all $\lambda\in (0,\pi)$
\begin{equation}
 \lim_{\substack{ t\to  \infty\\ t\to 0 } }t^{1/2 \pm ik}({\bf R} (\lambda+ i0) f)(t)=  \frac{  i  }{\lambda \sqrt{\pi^2 - \lambda^2}} \int_{0}^\infty s^{-1/2\pm  ik } f(s) ds.
\label{eq:Li3x}\end{equation}
 \end{proposition}
 
   \medskip

{\bf 3.3.}
 Next,  we describe singularities of the resolvent ${\bf R}(z) $  as $z$ approaches the edge point $z=\pi$. Let  ${\bf k} (z)$ be function   \e{eq:kk}. According to \e{eq:Kk1} we have  ${\bf k} (\pi) =2i$. Observe that
\begin{equation}
\theta (z):= 2+ i {\bf k} (z)=   \pi ^{-2}\sqrt{z^2-\pi^2}  + O(|z-\pi|) 
\label{eq:L1sy}\end{equation}
as $z \to \pi$.
In terms of $\theta (z)$, function  \e{eq:B4xx} can be written as
\begin{equation}
  \rho (u; z)=\frac{u^{\theta(z)} }{1-u^{ 2}}+ \frac{u^{-\theta(z)}}{1-u^{-2}}.
\label{eq:L1sz}\end{equation}
It follows   that for all fixed $u\in{\Bbb R}_{+}$ and $z \to \pi$
\[
\rho  (u; z)=  \sum_{n=0}^\infty \frac{ \theta(z)^n}{n!} \sigma_{n}(u)
\]
where $\sigma_{n}(u)= \ln^n u$ for even $n$ and
\begin{equation}
\sigma_{n}(u)=  \frac{1+ u^2 }{1-u^2}\ln^n u
\label{eq:L1s}\end{equation}
 for odd $n$. Obviously,
\[
|\sigma_{n}(u)|\leq C_{n} (1+|\ln u|)^n
\]
for all $n$. This yields     the following result.

  \begin{proposition}\label{pi}
  Let $t,s\in{\Bbb R}_{+}$ be fixed.
  As $z\to\pi$, kernel \e{eq:B5y} admits the expansion in the asymptotic series
  \begin{equation}
{\bf a} (t,s; z)=  \frac{ 1}{\sqrt{z^2 - \pi^2}} (ts)^{-1/2} \sum_{n=0}^\infty \frac{\theta(z)^n}{n!} \sigma_{n}(t/s).
\label{eq:L1szz}\end{equation}
 \end{proposition}

 Putting together  Theorem~\ref{res}  and Proposition~\ref{pi}, we see  that the integral kernel of the resolvent ${\bf R}(z) $   has the singularity 
$$ -\pi^{-1}  (z^2-\pi^2)^{-1/2} (ts)^{-1/2}$$
 as $z\to \pi$.
 Thus the   Carleman operator $\bf C$ has a resonance at the point $z=\pi$. Note that $\psi_{0}(t)=t^{-1/2}$ is the ``eigenfunction" of the operator $\bf C$ corresponding to the spectral point $ \pi$.  It satisfies the equation ${\bf C}\psi_{0}=\pi \psi_{0}$ and ``almost belongs" to $L^2 ({\Bbb R}_{+})$.

\medskip

{\bf 3.4.}
Finally,  we study the resolvent ${\bf R}(z) $ or, equivalently, the operator-valued function  ${\bf A}(z)$ as $z$ approaches the edge point   $z=0$. For simplicity, we suppose that $z=\lambda<0$. We proceed from representation \e{eq:RR} for the kernel of the operator ${\bf A}(\lambda)$. By definition \e{eq:B6k}, we have the asymptotic relation
\[
k(|\lambda|)= -\pi^{-1} \ln |\lambda| +\pi^{-1} \ln (\pi+\sqrt{\pi^2-\lambda^2})=
-\pi^{-1} \ln |\lambda| +\pi^{-1} \ln (2\pi)+ O(\lambda^2)
\]
as $\lambda\to - 0$.
 It follows that for every fixed $u\in{\Bbb R}_{+}$  
\[
\sin\big(k(|\lambda|)\ln u\big)= -\sin \big(\pi^{-1} \ln ( |\lambda|/  (2\pi))\ln u\big)+ O(\lambda^2).
\]
  Thus we obtain the following result.

\begin{proposition}\label{zero}
Let $t,s\in{\Bbb R}_{+}$, $t\neq s$, be fixed.
  As $\lambda\to -  0$,  the kernel of the operator ${\bf A}(\lambda)$ obeys   the asymptotic relation
  \begin{equation}
{\bf a} (t,s; \lambda)= \frac{ 2  }{\sqrt{ \pi^2-\lambda^2}}
\frac {\sqrt{ ts}}{t^2-s^2 }
 \sin \big(\pi^{-1} \ln ( |\lambda|/  (2\pi))\ln (t/s)\big)+ O(\lambda^2)   .
\label{eq:RRz}\end{equation}
 \end{proposition}

 Now formula \e{eq:B5} shows that the singularity of the resolvent at the point $z=0$ consists of the singular denominator $z^{-1}$ and of the oscillating term. Using formulas 
 \e{eq:B4xxx} and \e{eq:B5yy}   we can obtain similar results for  $z=\lambda\pm i0\to 0$ along the continuous spectrum.

  \section{Time-dependent evolution}  

\medskip

{\bf 4.1.}
Here we study the unitary group $\exp(-i {\bf C}T)$. It follows from relation \e{eq:C4} and Lemma~\ref{mellin} that
\[
\exp(-i {\bf C}T)=I + B(T)
\]
where 
\[
(B(T) f) (t)=\int_{0}^\infty (ts)^{-1/2} b(t/s; T) f(s)ds,
\]
\[
b(u; T)= (2\pi)^{-1} \int_{-\infty}^\infty u^{ik} (e^{-i \lambda (k) T} -1) dk
\]
and the function $ \lambda (k)$ is defined by formula \e{eq:C6}.
Apparently  the integral   $b(u; T)$ cannot be expressed  in terms of standard functions.

\medskip

{\bf 4.2.}
However it is possible to find explicitly the asymptotics of   $\exp(-i {\bf C}T)$ as $T\to \pm \infty$.  According to formula \e{eq:C4}  we have
\begin{equation}
(e^{-i {\bf C}T} f)(t)= (2\pi)^{-1/2} t^{-1/2 } \int_{-\infty}^\infty  e^{i(k \ln t-  \lambda(k)T)} \ti{f}(k)dk, \q \ti{f} =Mf.
\label{eq:StPh}\end{equation}
  Let us apply the stationary phase method to integral \e{eq:StPh}. Stationary points $k$  are determined by the equation
\begin{equation}
\lambda' (k)= \frac{\ln t}{T }  =: - \frac{\pi^2}{2 }  \tau.
\label{eq:StPh1}\end{equation}
Obviously, the function 
\begin{equation}
\lambda' (k)=- \pi^2 \frac{\sinh(\pi k)}{1+\sinh^2(\pi k)} 
\label{eq:StPh2}\end{equation}
is odd, it is negative for $k>0$ and $\lambda' (k)\to 0$ as $k\to \infty$. It has the minimum $- \pi^2/2$ at the point $k_{0}= \pi^{-1} \ln ( \sqrt{2}+1)$ and $\lambda'' (k) < 0$ for $k\in [0, k_{0})$, $\lambda'' (k) > 0$ for $k > k_{0}$. It follows that
equation \e{eq:StPh1} does not have solutions if $|\tau|> 1$,  it has two positive
solutions $k_{1}(\tau) < k_{0} <k_{2}(\tau)$ for $\tau\in (0, 1)$ and  it has two negative solutions  $k_{2}(\tau) < - k_{0} < k_{1}(\tau)$ for $\tau\in (- 1,0)$. Clearly, $k_{j}(-\tau) = - k_{j}(\tau)$. Let us introduce the short-hand notation
\[
 \sigma_{j}  ( \tau) =1+(-1)^j \sqrt{1-\tau^2}.
\]
An easy calculation shows that
\begin{equation}
\sinh(\pi k_{j}(\tau))=   \sigma_{j}  ( \tau)/ \tau
\label{eq:tks}\end{equation}
and
\begin{equation}
k_{j}(\tau)=  \pi^{-1} \sgn \tau \ln\big(|\tau|^{-1}(\sigma_{j} ( \tau)+\sqrt{2\sigma_{j} ( \tau)}) \big).
\label{eq:tk}\end{equation}

Let ${\cal M}$ be the set ${\Bbb R } $ with the points $0$, $k_{0}$ and $- k_{0}$ removed.   Applying the stationary phase method to integral \e{eq:StPh} where $\ti{f}\in C_{0}^\infty ({\cal M})$, we see that, for  
$e^{-\pi^2 |T |/2} < t <  e^{\pi^2 |T |/2}$, 
\begin{equation}
(e^{-i {\bf C}T} f)(t)=  |T|^{-1/2}  t^{-1/2 }   \sum_{j=1}^2 \d_{j} e^{- i \omega_{j}(\tau) T} |\lambda'' (k_{j}(\tau))|^{-1/2}  \ti{f}(k_{j}(\tau))  + O(|T|^{-1}).
\label{eq:StPh3}\end{equation}
Here  $\d_1=e^{ i    (\sgn T)\pi/4 }$, $\d_2=e^{- i  (\sgn T)\pi/4 }$ and 
\begin{equation}
\omega_{j} (\tau)=\pi^2 k_{j}(\tau) \tau/2 + \lambda(k_{j}(\tau)).
\label{eq:StPh5}\end{equation}
Note that $\tau\in (- 1, 0)$ for $t\in (1, e^{\pi^2 |T |/2}) $ and $\tau\in (0, 1)$ for $t\in ( e^{-\pi^2 |T |/2} ,1)$. If $|\tau|\geq 1$, that is, $t \geq  e^{\pi^2 |T |/2}$ or $t \leq  e^{-\pi^2 |T |/2}$, then integrating by parts, we see that integral \e{eq:StPh} decays faster than any power of $(|\ln t| + |T|)^{-1}$. Using formulas \e{eq:tks} and \e{eq:tk}, it is easy to calculate
 \begin{equation}
\lambda(k_{j}(\tau))= 2  \pi |\tau | \frac{\sigma_{j}(\tau)+ \sqrt{2\sigma_{j}(\tau)}  }
{ (\sigma_{j}(\tau)+ \sqrt{2\sigma_{j}(\tau)})^2+\tau^2}
\label{eq:tk3}\end{equation}
and
\begin{equation}
\lambda'' (k_{j}(\tau))=(-1)^j   \pi^3 |\tau| \sqrt{\frac{ 1-\tau^2} {2\sigma_{j}(\tau)}}.
\label{eq:tk2}\end{equation}

Let us formulate the result obtained. 
 
  \begin{lemma}\label{Time1}
 Suppose that $\ti{f}\in C_{0}^\infty ({\cal M})$. If $t \geq  e^{\pi^2 |T |/2}$ or $t \leq  e^{-\pi^2 |T |/2}$, then, for all $N$, we have the estimates 
  \[
|(e^{-i {\bf C}T} f)(t)|\leq C_{N}(|\ln t| + |T|)^{-N}.
\]
 If $t \in (  e^{-\pi^2 |T |/2},  e^{\pi^2 |T |/2})$, then the asymptotics of
 the function $(e^{-i {\bf C}T} f)(t)$ is given by formula \e{eq:StPh3} where $\tau$ and $k_{j}(\tau)$ are defined by relations \e{eq:StPh1} and \e{eq:tk}. The functions  
 $\omega_{j}(\tau)$, $\lambda (k_{j}(\tau))$ and $\lambda'' (k_{j}(\tau))$ are determined by equalities \e{eq:StPh5}, \e{eq:tk3} and \e{eq:tk2}, respectively.
  \end{lemma}

\medskip

{\bf 4.3.}
Let us now set
\begin{equation}
(U_{j} (T) f)(t)= \chi_{T} (t) |T|^{-1/2}  t^{-1/2 }    \d_{j} e^{ - i \omega_{j}(\tau) T} |\lambda'' (k_{j}(\tau))|^{-1/2}  \ti{f}(k_{j}(\tau))
\label{eq:StPh4}\end{equation}
where $\chi_{T} $ is the characteristic function of the interval $(e^{-\pi^2 |T |/2}, e^{\pi^2 |T |/2})$ and $\tau$ is related to $t$ and $T$ by equality \e{eq:StPh1}.    Let us  calculate
\begin{align}
\| U_{j} (T) f \|^2=& |T|^{-1 } \int_{e^{-\pi^2 |T |/2}}^{e^{\pi^2 |T |/2}} t^{-1  }|\lambda'' (k_{j}(\tau))|^{-1} | \ti{f}(k_{j}(\tau))|^2 dt 
\nonumber\\
= &\frac{\pi^2}{2} \int_{ - 1}^{ 1}  |\lambda'' (k_{j}(\tau))|^{-1} | \ti{f}(k_{j}(\tau))|^2 d\tau.
\label{eq:StPh6}\end{align}
According to \e{eq:StPh1}  we have $\lambda'(k_{j}(\tau))=-\pi^2\tau/2$.
Differentiating this equation we see that $\lambda''(k_{j}(\tau))k_{j}'(\tau)=-\pi^2/2$ and hence making in \e{eq:StPh6} the change of variables $k=k_{j}(\tau)$, we find that
\begin{equation}
\| U_{j} (T) f \|^2=  \int_{I_{j} }  | \ti{f}(k )|^2 dk
\label{eq:StPh8}\end{equation}
where $I_{1}=(-k_{0}, k_{0})$ and $I_2=(-\infty,- k_{0}) \cup (k_{0}, \infty)$. It follows that 
\begin{equation}
\| U_1(T) f \|^2 + \| U_2 (T) f \|^2=  \|   f \|^2 
\label{eq:StPh7}\end{equation}
for all $T$. In particular, the operators $U_{j} (T)$ are bounded uniformly in $T$.

Lemma~\ref{Time1} implies that 
 \begin{equation}
\lim_{|T|\to \infty} \| \big(e^{-i {\bf C}T}  -U_1(T)  -U_2(T) \big)f \|=0
\label{eq:StPhX}\end{equation}
 for    $\ti{f}\in C_{0}^\infty ({\cal M})$. Of course this relation extends to all $f\in L^2 ({\Bbb R}_{+})$ which yields   the following assertion.
 
  \begin{theorem}\label{Time}
 Define the operators $U_{j} (T)$, $j=1,2$, by formula  \e{eq:StPh4}. Then for all $f\in L^2 ({\Bbb R}_{+})$ relation \e{eq:StPhX} holds.
  \end{theorem}
  
  Note that the operators  $U  (T)=U_1 (T)+U_2 (T)$ are not unitary but $\|U(T)f\|\to 
 \| f\|$ for all $f\in L^2 ({\Bbb R}_{+})$  as $|T|\to \infty$.  This fact follows from Theorem~\ref{Time}. It is equivalent to the relation
   \begin{equation}
\lim_{| T |\to \infty} \Re \,( U_1(T) f  , U_2(T) f ) =  0
\label{eq:StPh9}\end{equation}
which can be verified directly with the help of representation  \e{eq:StPh4}.

\section{Scattering theory}  

{\bf 5.1.}
Our   goal now  is to develop   scattering theory for perturbations of the Carleman operator $H_{0} = {\bf C}$ by self-adjoint operators $V$ satisfying the condition
\begin{equation}
Q^{\alpha} V Q^{\alpha}\in {\goth S}_{\infty}
\label{eq:H4c}\end{equation}
for some $\alpha>1/2$. Here the operator $Q$ is   defined by formula \e{eq:QQ}, and ${\goth S}_{\infty}$ is the class of compact operators. For example, condition \e{eq:H4c} holds if $V$ is an integral operator
\begin{equation}
(V f)(t)=\int_{0}^\infty {\bf v}  (t,s)  f(s) ds,
\label{eq:v1}\end{equation}
with kernel ${\bf v}  (t,s)$ such that  
\begin{equation}
\int_0 ^\infty \int_{0}^\infty | {\bf v}  (t,s) |^2 \la \ln t\ra^{2\alpha}\la \ln s\ra^{2\alpha}   dt ds<\infty. 
\label{eq:H4}\end{equation}
Then the operator $Q^{\alpha} V Q^{\alpha}$ belongs to the Hilbert-Schmidt class.

 In particular, if ${\bf v}  (t,s) = v (t+s)$, then 
  $V$ is a Hankel operator acting by formula \e{eq:H1}. Since
\[
\int_0 ^t   \la \ln s\ra^{2\alpha}    \la \ln (t-s) \ra^{2\alpha}  ds\leq C   \la \ln t\ra^{4\alpha} t,
\]
  condition \e{eq:H4}  is satisfied for a Hankel operator $V$ if
\begin{equation}
\int_0 ^\infty  | v (t)|^2 \la \ln t\ra^{4\alpha}  t  dt  <\infty . 
\label{eq:H4h}\end{equation}

Main results of scattering theory for the pair of  the operators $ H_{0}={\bf C}$
 and $H=H_{0}+V$ are collected in the following assertion.

    \begin{theorem}\label{ScTh}
    Let  assumption \e{eq:H4c} hold for some $\alpha>1/2$. Then:
    
    $1^0$ The strong limits
        \begin{equation}
\slim_{T \to\pm\infty} e^{iH T} e^{- iH_{0}T} =: W_{\pm} 
\label{eq:W1}\end{equation}
known as the wave operators exist.

    $2^0$
 The wave operators enjoy the intertwining property  $H W_{\pm}  =W_{\pm} H_{0}$ and are isometric.

 $3^0$ The  wave operators are  complete:
  \begin{equation}
\Ran  W_{\pm}   =     {\cal H}^{(\ac)}
\label{eq:W6}\end{equation}
where ${\cal H}^{(\ac)}$ is the absolutely continuous subspace  of the operator $H$.
     \end{theorem}
  
\begin{corollary}\label{scth}
 The absolutely continuous spectrum  of the 
operator $H  $ has multiplicity $2$ and coincides with the interval $[0,\pi]$.
\end{corollary}

A proof of Theorem~\ref{ScTh} can be obtained by means of the stationary scattering theory. Denote $R_{0} (z)=(H_{0}-z I )^{-1}$  (of course, $R_{0} (z)= {\bf R}(z)$ in the notation of Section~2),    $R  (z)=(H -z I)^{-1}$ and recall the resolvent identity 
 \begin{equation}
 R (z)   - R_{0}(z) =-    R (z)  VR_{0}(z) =-    R_{0} (z)  VR (z).
\label{eq:RESx}\end{equation}
Set $G_{0}(z) = Q^{-\beta} R_{0}(z)    Q^{-\beta}$ and $G (z) = Q^{-\beta} R (z)    Q^{-\beta}$. It follows from  \e{eq:RESx} that
 \begin{equation}
\big( I+ G_{0}(z) K   \big)         G (z)     = G_{0}(z)     
\label{eq:RES}\end{equation}
where $K=Q^{\beta} V    Q^{\beta}\in {\goth S}_{\infty}$ for all $\beta\leq \alpha$. Considering \e{eq:RES} as an equation for the operator-valued function $ G (z) $ and using Proposition~\ref{cont} (the limiting absorption principle for the operator $H_{0}$), one can prove (see, e.g., Theorem~7.3 of  Chapter 4 of \cite{Yafaev}) a similar statement   for the operator $H$. To formulate the precise result, let us introduce the set       ${\cal N}$   of $\lambda\in (0,\pi)$ where at least one of the equations 
     \[
  f+   G_{0}(\lambda\pm i0) Kf=0
  \]
    has a nontrivial solution.

 \begin{theorem}\label{LAP}
     Let  assumption \e{eq:H4c} hold for some $\alpha >1/2 $.  Then the set ${\cal N} $   is closed and has the Lebesgue measure zero.   The  singular   spectrum of the operator $H$ is contained in the set $\cal N$. The operator-valued function $G(z)= Q^{-\beta} R(z) Q^{-\beta} $ where $\beta>1/2$ is H\"older continuous with any exponent $\gamma< \min\{\alpha,\beta\}-1/2$ $($and $\gamma\leq 1)$  in $z$ if $\pm\Im z\geq 0$ and $\Re z\in (0,\pi)\setminus {\cal N}$. 
\end{theorem}

Theorem~\ref{ScTh} can be deduced from Theorem~\ref{LAP} 
 (see, e.g., Theorems~6.4 and 6.5 of  Chapter 4 of \cite{Yafaev}, for details). 
 
If $\beta >1$, then the  operator-valued function $G_{0}(z)  $ is H\"older continuous with an exponent $\gamma >1/2$. This allows one to show  (see, e.g., Theorems~7.9 and 7.10 of  Chapter 4 of \cite{Yafaev}) that  the  singular set ${\cal N}$ coincides with the point spectrum $\spec_{p} H$ of $H$ and leads to the following result.
  
 \begin{theorem}\label{spectral}
     Let  assumption \e{eq:H4c} hold for some $\alpha >1 $.  Then  the operator $H=H_{0}+V$
  does not have the  singular continuous spectrum.
   Eigenvalues of   $H$ distinct from the points $0$ and $\pi$ have
  finite multiplicities and can   accumulate to these points only. 
  \end{theorem}

  \medskip

{\bf 5.2.}
Standard formulas of   stationary scattering theory allow us to obtain an expansion over eigenfunctions of the continuous spectrum of the operator $H$.  For $k>0$, we set 
$$
\psi_{1}^{(0)}(t; k)=t^{-1/2+ ik}, \q \psi_2^{(0)} (t ;k)= t^{-1/2- ik}
$$ 
  and define eigenfunctions  $\psi_{j}^{(\pm)}(t, k)$ of the operator $H$ by the formula
  \begin{equation}
 \psi_{j}^{(\pm)}(k)= \psi_{j}^{(0)}(k) -R(\lambda (k) \mp i 0) V  \psi_{j}^{(0)}(k),\q j=1,2,
 \q \lambda(k)\not\in \spec_{p} H,
\label{eq:E}\end{equation}
where $\lambda= \lambda(k)\in (0,\pi)$ and $k \in{\Bbb R}_{+}$ are related by formula \e{eq:C6}.    It follows from Theorems~\ref{LAP} and \ref{spectral} that 
 \begin{equation}
 \la \ln t\ra^{- \beta} \psi_{j} (t, k)\in L^2({\Bbb R}_{+}), \q \beta >1/2,  
 \label{eq:Ey}\end{equation}
 and these functions depend continuously (actually, H\"older continuously) on $k>0$ in the norm of this space. Obviously, we have $H\psi_{j}^{(\pm)}(k)= \lambda (k) \psi_{j}^{(\pm)}(k)$.


Set 
 \begin{equation}
  (\Psi_0 f)(k)=   \big( (M f)(k),(M f)(-k) \big)^\top
\label{eq:E1a}\end{equation}
where  $M$ is  the Mellin transform  defined by formula \e{eq:C3}. Clearly, the mapping $\Psi_0 : L^2({\Bbb R}_{+})\to 
L^2({\Bbb R}_{+}; {\Bbb C}^2)$ is unitary and according to \e{eq:C4} it diagonalizes the operator $H_{0}$:
 \begin{equation}
(\Psi_{0} H_{0} f)(k) = \lambda (k) ( \Psi_{0} f)(k), \q \forall f\in L^2({\Bbb R}_{+}).
\label{eq:E1ax}\end{equation}
Let us now construct a diagonalization of the operator $H$.

 \begin{theorem}\label{Psi}
 Define operators $\Psi_\pm $  by the relation
 \begin{equation}
  (\Psi_\pm f)(k)= (2\pi)^{-1/2} \Big( \int_{0}
^\infty \overline{ \psi_1^{(\pm)}(t ; k)} f(t) dt, \int_{0}
^\infty \overline{ \psi_2^{(\pm)}(t ; k)} f(t) dt \Big)^\top
\label{eq:E1}\end{equation}
  on functions $f$ such that $\la \ln t\ra^{\beta} f \in L^2({\Bbb R}_{+})$ for some $ \beta >1/2$. Then $\Psi_\pm f\in L^2({\Bbb R}_{+}; {\Bbb C}^2)$ and the operators $\Psi_\pm$ extend to bounded operators $\Psi_\pm: L^2({\Bbb R}_{+})\to L^2({\Bbb R}_{+}; {\Bbb C}^2)$. They satisfy the relations
   \begin{equation}
 \Psi_\pm ^*  \Psi_\pm =I- E^{({\rm p})}, \q  \Psi_\pm  \Psi_\pm^*=I
\label{eq:E2}\end{equation}
where $E^{({\rm p})}$ is the orthogonal projection on the subspace ${\cal H}^{({\rm p})}$ spanned by eigenvectors of the operator $H$. The operators $\Psi_\pm $ diagonalize $H$, that is,
\begin{equation}
( \Psi_\pm H f)(k) = \lambda(k) ( \Psi_\pm  f)(k), \q \forall f\in L^2({\Bbb R}_{+}),
\label{eq:E3}\end{equation}
and   are related to the wave operators $W_\pm $ by the equality
 \[
W_{\pm}= \Psi_\pm^* \Psi_{0}.  
\]

\end{theorem}

 \medskip

{\bf 5.3.}
In terms of the wave operators \e{eq:W1}, the scattering operator is defined by the formula
$
{\bf S}=W_{+}^* W_{-}.
$
According to Theorem~\ref{ScTh} the scattering operator commutes with $H_{0}$, that is, ${\bf S} H_{0}=H_{0} {\bf S}$,  and it is unitary.

To  define the scattering matrix, we use the diagonalization \e{eq:E1ax}  of the operator $H_{0}$.    Since the operators ${\bf S} $ and $H_{0}$ commute, we have
\[
( {  \Psi}_0 {\bf S}  f)(k)= S(k)( { \Psi}_0    f)(k), \q k>0,
\]
where the $2 \times 2$ matrix 
  \begin{equation}
S(k)= \begin{pmatrix}
 s_{11}(k) &  s_{12}(k)
 \\ 
 s_{21}(k) &  s_{22}(k)
 \end{pmatrix} 
\label{eq:SMM}\end{equation}
 is known as the scattering matrix. The matrices $S(k)$ are unitary for all $k\in {\Bbb R}_{+}$ because the operator ${\bf S} $ is unitary.  Observe that we parametrize the scattering matrix $S(k)$      by the quasi-momentum $k$.
  
The scattering matrix can be expressed  in terms of boundary values of the resolvent $R(z)$.   Set   
  \begin{equation}
\Gamma_{0}(k)f=  \sqrt{\gamma(k))} \big( (M f)(k ),(M f)(-k ) \big)^\top
\label{eq:E1ab}\end{equation}
where
\begin{equation}
 \gamma(k)= |\lambda' (k) |^{-1}= \frac{\cosh^2(\pi k)}{ \pi^2 \sinh(\pi k)}.
\label{eq:smc}\end{equation}
Then
\[
\int_{0}^\pi \| \Gamma_{0}(k (\lambda))f\|_{{\Bbb C}^2}^2 d\lambda=
\int_{-\infty}^\infty |  (M f)(k)|^2 d k= \| f\|^2.
\]
Using, e.g., Theorem~5.5.3 and Proposition~7.4.1   of \cite{Yafaev}, we can state the following result. 

\begin{proposition}\label{Sst}
    Let assumption \e{eq:H4c} hold  for some $\alpha>1$. Then the scattering matrix  admits the representation
 \begin{equation}
S(k) = I-2\pi i \Gamma_{0}(k) \big(V-V R(\lambda (k)+i0) V \big)\Gamma_{0}^*(k), \q \lambda (k)\not\in\spec_{p} H .
\label{eq:SS}\end{equation}
\end{proposition}
 
Observe that the right-hand side of \e{eq:SS} can be written as a combination of bounded operators.

  Formulas \e{eq:E1} -- \e{eq:E3} mean that the functions $\psi_{j}^{(+)}(t;k)$ and $\psi_{j}^{(-)}(t;k)$ where $j=1,2$ give two ``bases" in the ``eigenspace" of $H$ corresponding to the ``eigenvalue"  $\lambda(k)$. Since the absolutely continuous spectrum of the operator $H$ has multiplicity $2$, it is natural to expect that the functions $\psi_{j}^{(-)}(t;k)$, $j=1,2$, are linear combinations of the functions $\psi_1^{(+)}(t;k)$ and $\psi_2^{(+)}(t;k)$ (and vice versa). It turns out that the link between these two ``bases" is given by the elements of scattering matrix \e{eq:SMM}. The proof of the following assertion is very similar to the corresponding result of \cite{LNM}.
  
  \begin{proposition}\label{pm}
    Let assumption \e{eq:H4c} hold  for some $\alpha>1$. Then
      \begin{equation}
      \begin{split}
 \psi_{1}^{(-)} (t;k)&= s_{11}(k)\psi_{1}^{(+)} (t;k)+ s_{21}(k)\psi_{2}^{(+)} (t;k)
 \\
  \psi_{2}^{(-)} (t;k)&= s_{12}(k)\psi_{1}^{(+)} (t;k)+ s_{22}(k)\psi_{2}^{(+)} (t;k).
  \end{split}
\label{eq:pm}\end{equation}
\end{proposition}

\begin{pf}
It can be deduced from the resolvent identity \e{eq:RESx} and formula \e{eq:SS} for the scattering matrix that  
  \begin{equation}
    \big(I- R(\lambda (k)+i0) V \big)\Gamma_{0}^*(k)=    \big(I- R(\lambda (k)-i0) V \big)\Gamma_{0}^*(k)   S(k).
\label{eq:pm1}\end{equation}
Let $b=(b_{1}, b_{2})^\top \in {\Bbb C}^2$. It follows from  \e{eq:E1ab} that
\[
(\Gamma_{0}^*(k) b)(t)=\sqrt{\gamma(k)} (2\pi)^{-1/2} (b_{1}t^{-1/2+ik}
+ b_2 t^{-1/2- ik})
\]
and hence, by definition \e{eq:E},
 \begin{equation}
((I- R(\lambda (k)\pm i0) V ) \Gamma_{0}^*(k) b)(t)=\sqrt{\gamma(k)} (2\pi)^{-1/2} (b_{1} \psi_{1}^{(\mp)}(t;k)
+ b_2 \psi_{2}^{(\mp)} (t;k)).
\label{eq:pm2}\end{equation}
Thus relation \e{eq:pm1} implies that
 \begin{multline*}
b_{1} \psi_{1}^{(-)}(t;k)
+ b_2 \psi_{2}^{(-)} (t;k)
\\
=(s_{11}(k)b_{1}+ s_{12}(k)b_{2}) \psi_{1}^{(+)}(t;k)
+ (s_{21}(k)b_{1}+ s_{22}(k)b_{2})  \psi_{2}^{(+)} (t;k).
 \end{multline*}
Comparing here the coefficients at $b_{1}$ and $b_{2}$, we arrive at formulas 
\e{eq:pm}.
\end{pf}

Using the unitarity of $S(k)$, we can rewrite formulas \e{eq:pm} as
\[
      \begin{split}
 \psi_{1}^{(+)} (t;k)&=\overline{s_{11}(k)}\psi_{1}^{(-)} (t;k)+ \overline{s_{12}(k)}\psi_{2}^{(-)} (t;k)
 \\
  \psi_{2}^{(+)} (t;k)&= \overline{s_{21}(k)}\psi_{1}^{(-)} (t;k)+ \overline{s_{22}(k)}\psi_{2}^{(-)} (t;k).
  \end{split}
\]

 \medskip

{\bf 5.4.}
The representation  \e{eq:SS} of the scattering matrix  $S(k)$ can be rewritten in terms of eigenfunctions $\psi_{j}(t ; k): =\psi_{j}^{(-)}(t ; k)$, $j=1,2$, of the operator $H$.  

\begin{proposition}\label{SS}
   Let assumption \e{eq:H4c} hold  for some $\alpha>1$. Then the elements of  the scattering matrix \e{eq:SMM}   are given by the formulas 
  \begin{align}
  s_{11}(k)&=1-i \gamma(k)\int_0 ^\infty t^{-1/2- ik} (V \psi_{1}(k) )(t)dt ,
\label{eq:sm11}\\
  s_{12}(k)&=- i \gamma(k)\int_0^\infty t^{-1/2+ik} (V \psi_{1}(k) )(t)dt ,
\label{eq:sm12}\\
  s_{21}(k)&=-i\gamma(k)\int_0^\infty t^{-1/2-ik} (V \psi_{2}(k) )(t)dt ,
\label{eq:sm21}\\
  s_{22}(k)&=1-i\gamma(k)\int_0^\infty t^{-1/2+ik} (V \psi_{2}(k) )(t)dt ,
\label{eq:sm22}\end{align}
where the coefficient $ \gamma(k)$ is defined by formula \e{eq:smc}.
\end{proposition}

\begin{pf} 
 Let again $b=(b_{1}, b_{2})^\top$.  It follows from definition \e{eq:E1ab}  and representation \e{eq:pm2} that
  \begin{align*}
2\pi i \Gamma_{0} (k) V (I- R(\lambda (k)+i0) V ) \Gamma_{0}^*(k) b 
=i \gamma(k) \big(F_{+}(k) b ,  F_{-}(k) b  \big)^\top
 \end{align*}
 where we have used the notation
  \begin{equation}
F_{\pm}(k) b = \int_0^\infty t^{-1/2\mp ik} (V (b_{1} \psi_{1}(k)
+ b_2 \psi_{2}(k)) )(t)dt.
\label{eq:SF} \end{equation}
Now representation  \e{eq:SS} implies  that
\[
S(k)b=b-i \gamma(k) \big(F_{+}(k)b ,  F_{-}(k) b  \big)^\top.
\]
In view of \e{eq:SF} this formula for the matrix $S(k)$ is equivalent to formulas \e{eq:sm11} -- \e{eq:sm22} for its elements.
\end{pf}

Note that the matrices $S(k)$   depend continuously (actually, H\"older continuously) on $k>0$ (away from the point spectrum of $H$).

Our next goal is to find  asymptotics of the functions $\psi_{j} (t; k)$ as $t\to \infty$ and as $t\to 0$. We proceed from  the Lippmann-Schwinger equation
  \begin{equation}
 \psi_{j} (k)= \psi_{j}^{(0)}(k) -R_{0}(\lambda (k) + i 0) V  \psi_{j} (k),\q j=1,2,
\label{eq:Li}\end{equation}
for the functions $\psi_{j} (k)$. Recall that this equation is an immediate consequence of the resolvent identity \e{eq:RESx} and the definition  \e{eq:E} of $\psi_{j} (k)$.

  \begin{proposition}\label{S}
    Let assumption \e{eq:H4c} hold  for some $\alpha>1$, and let
     \begin{equation}
  (V Q^\alpha g)(t)= o(t^{-1/2}), \q \forall g\in L^2 ({\Bbb R}_{+}), 
\label{eq:V1}\end{equation}
as $t \to \infty$ and as $t\to 0$
for some $\alpha>1/2$.  Then
 \begin{equation}
  \begin{split}
\psi_{1}(t; k)&= s_{11}(k)t^{-1/2+ik}+ o(t^{-1/2}), \q t\to 0,
 \\
\psi_{1}(t;k)&=t^{-1/2+ik} +  s_{12}(k) t^{-1/2- ik}+ o(t^{-1/2}), \q t\to \infty,
 \end{split}
\label{eq:as1}\end{equation}
and
 \begin{equation}
  \begin{split}
\psi_{2}(t;k)&=t^{-1/2- ik} + s_{21}(k)t^{-1/2+ik}+ o(t^{-1/2}), \q t\to 0,
 \\
\psi_{2}(t; k)&=   s_{22}(k)t^{-1/2-ik}+ o(t^{-1/2}), \q t \to \infty,
 \end{split}
\label{eq:as2}\end{equation}
where the  asymptotic coefficients are the elements of scattering matrix \e{eq:SMM}. 
\end{proposition}

 \begin{pf}
  Set $f_{j}(t)=(V   \psi_{j} (k))(t) $. In view of of equation \e{eq:Li}  we only have to find  asymptotics of the functions $(R_{0}(\lambda (k) + i 0)f_{j})(t)$ as $t \to \infty$ and as $t\to 0$. Combining assumption  \e{eq:H4c} for $\alpha>1/2$ and inclusion  \e{eq:Ey}, we see that     the functions $f_{j}(t)  $ obey  condition \e{eq:Li2x}. Moreover, $f_{j}(t)=  o (t^{-1/2})$ if condition \e{eq:V1}  is satisfied.
 Therefore  applying     Proposition~\ref{asympt1x} to functions $f_{j}$, we see that
 \[
 (R_{0}(\lambda (k)+i0) V \psi_{j} (k))(t)=i \gamma(k)t^{-1/2 +ik}
 \int_{0}^\infty s^{-1/2 - ik} (V\psi_{j} (k))(s)ds + o (t^{-1/2  })
 \]
as    $t\to 0$ and
 \[
 (R_{0}(\lambda (k)+i0) V \psi_{j} (k))(t)=i \gamma(k)t^{-1/2 -ik}
 \int_{0}^\infty s^{-1/2 + ik} (V\psi_{j} (k))(s)ds + o (t^{-1/2  })
 \]
as    $t\to \infty$.
 Substituting these asymptotic relations into the right-hand side of equation \e{eq:Li} for $\psi_{j}(t;k)$ and using equalities \e{eq:sm11} -- \e{eq:sm22}, we get  formulas \e{eq:as1} and \e{eq:as2}.
 \end{pf} 
 
Note that for integral operators \e{eq:v1}   relation \e{eq:V1} is true if
 \begin{equation}
  \int_0 ^\infty |v(t,s)|^2 \la \ln s\ra^{2\alpha} ds=o(t^{-1})  
\label{eq:V1+}\end{equation} 
as $ t\to \infty $ and $ t\to 0$.
For  Hankel operators \e{eq:H1},  assumption \e{eq:H4h} implies  both conditions \e{eq:H4} and \e{eq:V1+}.  
 
 Similarly to the one-dimensional Schr\"odinger equation (see \cite{F}), asymptotic relations \e{eq:as1} and \e{eq:as2} can be interpreted in the following way. The solution $\psi_1(t; k)$ describes a wave propagating from $t=\infty$ to $ t= 0$. According to  \e{eq:as1} we observe the transmitted part  $s_{11}(k)t^{-1/2+ik}$ going to $0$ and the reflected part $s_{12}(k)t^{-1/2- ik}$ going back to $\infty$.  In the same way, the solution $\psi_2(t; k)$ describes a wave propagating from $t= 0$ to $ t= \infty$. It is natural to call   $s_{11}(k)$,  $s_{22}(k)$ the transmission coefficients and  to call   $s_{11}(k)$,  $s_{22}(k)$ the reflection coefficients. The squares $|s_{jl}(k)|^2$ give probabilities of the corresponding processes. As might be expected,  $|s_{11}(k)|^2 + |s_{12}(k)|^2=1$ and $|s_{21}(k)|^2 + |s_{22}(k)|^2=1$.

 \medskip

{\bf 5.5.}
We always suppose that the operators $V$ are self-adjoint which corresponds to the condition
\begin{equation}
{\bf v} (t,s)=\overline{{\bf v} (s,t)} 
\label{eq:SA}\end{equation}
for integral operators \e{eq:v1}.
Let   the complex conjugation ${\cal C}$ be defined  by   equality \e{eq:CC}.  
Assume additionally that the operator $V$  commutes with ${\cal C}$, that is, 
\begin{equation}
{\cal C}V=V{\cal C}
\label{eq:compl}\end{equation}
 and hence ${\cal C}H=H{\cal C}$.   In terms of kernels, it means that ${\bf v}(t,s)$ is a real function which in view of the self-adjointness condition \e{eq:SA} is equivalent to the equality
\begin{equation}
 {\bf v} (t,s)= {\bf v} (s,t).
\label{eq:smv}\end{equation}
For self-adjoint Hankel operators \e{eq:H1}, this condition   is always satisfied.

Under assumption \e{eq:compl} the wave operators \e{eq:W1} obey the identity ${\cal C}W_{\pm}=W_{\mp}{\cal C}$ and hence 
\begin{equation}
{\cal C} {\bf S}= {\bf S}^* {\cal C}. 
\label{eq:SAC}\end{equation}
Set ${\cal J}(b_{1}, b_{2})^\top = (b_2, b_1)^\top$. Since operator \e{eq:E1a} satisfies the identity ${\cal C}\Psi_{0}= {\cal J}\Psi_{0} {\cal C}$, it follows from \e{eq:SAC} that
\[
{\cal J}{\cal C} S(k)= S^*(k) {\cal J}{\cal C}
\]
for all $k>0$. It is easy to see that the last identity is equivalent to the equality
\begin{equation}
 s_{11}(k)= s_{22}(k), \q k>0.
 \label{eq:SAC1}\end{equation}
 
 Alternatively, identity \e{eq:SAC1} can be deduced from representation \e{eq:SS} if one takes into account that ${\cal C}R(z)= R(\bar{z}){\cal C}$.

Observe also that under assumption \e{eq:compl} eigenfunctions \e{eq:E} of the operator  $H$ are linked by the relations $\overline{\psi_{1}^{(+)} (t ; k)}=\psi_{2}^{(-)}(t ; k)$ and $\overline{\psi_{2}^{(+)} (t ; k)}=\psi_{1}^{(-)}(t ; k)$.

  \section{The discrete spectrum above the continuous spectrum}  

 Here we study the   spectrum of the perturbed Carleman operator $H=H_{0}+V$ lying above the point $\lambda=\pi$.  We prove   that under natural conditions on the kernel $v(t)$   it  consists of a finite number of eigenvalues. We also show that 
 for $V>0$ the operator $H$ has at least one eigenvalue larger than $\pi$.
In this section we suppose that the spectral parameter $\lambda \geq\pi$.

 \medskip
 
{\bf 6.1.} 
According to Proposition~\ref{pi}
  the nature of the singularity of the free resolvent $R_{0} (z)$ at the point $z=\pi$ is the same as that of the resolvent of the operator $ D^2 $ acting in the space $L^2 (\Bbb R)$ at the point $z=0$. So one can expect that the results on the spectrum of the perturbed Carleman operator $H=H_{0}+V$ above the point $ \pi$ are qualitatively similar to those on the negative spectrum of the Schr\"odinger operator $D^2+ {\sf V}(x)$. Here we verify this conjecture.

We proceed from   the Birman-Schwinger principle which in our case is formulated as follows.

 \begin{proposition}\label{B-S}
 Let $H_{0}$ be a bounded self-adjoint operator such that $H_{0}\leq\pi$.
Suppose that $V\geq 0$ and $V\in{\goth S}_{ \infty}$. Then   the total number $N(\lambda)$ of eigenvalues of the operator $H=H_{0}+V$ larger\footnote{Larger is strictly larger.} than $\lambda>\pi$ equals the total number   of eigenvalues of the operator $B(\lambda)=V^{1/2}(\lambda-H_{0})^{-1} V^{1/2}$  larger than $1$. 
 \end{proposition}

 \medskip
 
{\bf 6.2.}
 Representation \e{eq:B5} shows that 
 \begin{equation}
B(\lambda)=\lambda^{-1}\big(V+  V^{1/2}{\bf A}(\lambda) V^{1/2}\big).
\label{eq:M}\end{equation}
Recall that the integral kernel $ {\bf a} (t,s; \lambda) $ of the operator ${\bf A}(\lambda)$ is given by formula \e{eq:B5y} where $\rho(u;\lambda)$ is function \e{eq:L1sz}.   For $\lambda\geq\pi$,  function \e{eq:kk} equals 
\[
{\bf k}(\lambda)= \pi^{-1} i\arg (\pi- i\sqrt{\lambda^2-\pi^2}), 
\]
and  hence function \e{eq:L1sy} equals 
 \[
 \theta(\lambda)= \pi^{-1} \arctan (\pi^{-1}\sqrt{\lambda^2-\pi^2})\in [0,1/2).
 \]

Expansion  \e{eq:L1szz} shows that the kernel ${\bf a} (t,s; \lambda)$ is singular as $\lambda\to\pi$. Distinguishing the first (singular) term, we can rewrite formula \e{eq:B5y}   as
 \begin{equation}
{\bf a} (t,s; \lambda)= \frac{ 1}{\sqrt{\lambda^2 - \pi^2}} (ts)^{-1/2}+\ti{\bf a} (t,s; \lambda)
\label{eq:M1}\end{equation}
where
 \begin{equation}
 \ti{\bf a}(t,s; \lambda) = \frac{1}{\sqrt{\lambda^2 - \pi^2}} (ts)^{-1/2} \ti{\rho}(t/s ; \lambda)\big)
\label{eq:M2}\end{equation}
and
 \begin{equation}
   \ti{\rho} (u ; \lambda)= \rho (u ; \lambda)-1=
   \frac{u^{\theta(\lambda)}-1}{1-u^{2}} + \frac{u^{-\theta(\lambda)}-1}{1-u^{-2}} .
\label{eq:M2s}\end{equation}

 Let $\wt{\bf A} (\lambda)$ be the integral operator with kernel $ \ti{\bf a}  (t,s; \lambda)$. Let us show that it has the limit $ \wt{\bf A}  (\pi)=: \wt{A}$ as $\lambda\to \pi$.  We denote by $\| A\|_{2}$ the norm of  an operator $A\in {\goth S}_2$ in the  Hilbert-Schmidt class ${\goth S}_2$.

  \begin{lemma}\label{BS1x}
Let $\wt{A} $ be the   integral operator with kernel
 \begin{equation}
\ti{a}(t,s)=  \pi^{-2} (ts)^{-1/2} \sigma_{1} (t/s)
\label{eq:M4y}\end{equation}
where the function $\sigma_{1}$ is defined by   formula \e{eq:L1s}.
Then for all $\alpha>3/2$, we have
  \begin{equation}
\lim_{\lambda\to \pi} \| Q^{-\alpha}\big(\wt{\bf A}  (\lambda) - \wt{ A} \big)Q^{-\alpha}\|_2=0.
\label{eq:M41}\end{equation}
   \end{lemma}

 \begin{pf}
 We have to check that
  \begin{equation}
\lim_{\lambda\to \pi}\int_0^\infty \int_0^\infty \langle \ln t\rangle^{-2\alpha}
|\ti{\bf a} (t,s; \lambda) -\ti{ a} (t,s)|^2 \langle \ln s\rangle^{-2\alpha}dtds=0.
 \label{eq:L14}\end{equation}
  Since 
 \[
 \lim_{\lambda\to \pi} \frac{\ti{\rho}(u; \lambda)}{\sqrt{\lambda^2-\pi^2}}=
\pi^{-2}  \lim_{\lambda\to \pi}\big( \theta^{-1} (\lambda) \ti{\rho}(u; \lambda)\big) =  \pi^{-2} \sigma_{1} (u),
 \]
   the integrand in \e{eq:L14} tends to zero for all $t>0$ and $s>0$. We shall show that
  \begin{equation}
| \ti{\bf a}  (t,s; \lambda)  | \leq C (ts)^{-1/2} (1+ |\ln (t/s)| )
\label{eq:L17x}\end{equation}
where $C$ does not depend on $\lambda\geq \pi$. Then  the integrand in \e{eq:L14} is bounded by the function  
\[
C_{1} (ts)^{-1 }  \langle \ln t\rangle^{-2\alpha}\langle \ln s\rangle^{-2\alpha}
\big(\langle \ln t\rangle^{2}+\langle \ln s\rangle^{2}\big)
\]
 which belongs to $L^1 ({\Bbb R}^2)$ if $2\alpha>3$. 
 By the Lebesgue dominated convergence theorem,  this implies relation \e{eq:L14}.
 
 According to \e{eq:M2} for the proof of \e{eq:L17x}, we have to check that function \e{eq:M2s} satisfies the bound
 \begin{equation}
\ti{\rho}  (u; \lambda) \leq C \theta(\lambda) (1+|\ln u|), \q u \in{\Bbb R}_{+}.
\label{eq:L17}\end{equation}
Since $\ti{\rho} (u^{-1} ; \lambda)  = \ti{\rho}  (u; \lambda)$, 
it suffices to consider $u\geq 1$.  If $\theta \ln u\leq 1$, we use that
\[
\Big| \frac{u^{\pm\theta }-1}{1-u^{\pm 2}}\Big| \leq C \theta \frac{\ln u}{|1-u^{\pm 2}|}
\leq C_{1}   \theta (1+ \ln u),\q u\geq 1.
\]
 If $\theta \ln u\geq 1$ (and thus $   u\geq e^2$), then the absolute value of  the first  term in the right-hand side of \e{eq:M2s} is estimated by  $u ^{\theta  } ( u^{2  }-1)^{-1}$ which is bounded by a constant because $\theta < 1/2$. The   absolute value of the second term in the right-hand side of \e{eq:M2s}   is estimated by  $(1-u^{-2})^{-1}\leq(1-e^{-4})^{-1}$.  This concludes the proof of estimate \e{eq:L17}    and hence of relation \e{eq:L14}.
 \end{pf}

   Let us return to operator \e{eq:M}. Using \e{eq:M1}, we see that
 \begin{equation}
B(\lambda)=\frac{ 1}{\lambda\sqrt{\lambda^2 - \pi^2}}(\cdot, w)w +\wt{B}(\lambda)
\label{eq:M3}\end{equation}
where  $\psi_{0}(t) =t^{-1/2}$, $w=V^{1/2} \psi_{0}$ and  
 \[
\wt{B}(\lambda)=\lambda^{-1} \big(V +  V^{1/2} \wt{\bf A} (\lambda) V^{1/2} \big).
\]
Observe that under assumption \e{eq:H4c} the operator $V^{1/2}Q^{\alpha} $ is bounded and hence   $w\in L^2 ({\Bbb R}_{+})$ if $\alpha>1/2$. Let us   consider the operator $\wt{B}(\lambda)$. The following assertion is a direct consequence  of Lemma~\ref{BS1x}.

 \begin{lemma}\label{BS1}
Suppose that $V\geq 0$, $V\in {\goth S}_{2}$ and that inclusion \e{eq:H4c} is true for some $\alpha>3/2$. Then the operator $\wt{B}(\lambda)$ has the limit 
 \begin{equation}
\wt{B} (\pi)=\pi^{-1} \big( V +  V^{1/2} \wt{A}  V^{1/2}\big)
\label{eq:M4x}\end{equation}
 in the Hilbert-Schmidt norm as $\lambda\to \pi$. 
   \end{lemma}

 \medskip
 
{\bf 6.3.}
Now we return to the study of the discrete spectrum of the operator $H$.
   Let $\ti{n}(\lambda)$ (and $n (\lambda)$) be the total numbers   of eigenvalues of the operators $\wt{B}(\lambda) $ (resp. $B(\lambda)$)  larger than $1$. According  to relation \e{eq:M3} we have
   \[
    \ti{n} (\lambda)\leq n(\lambda)\leq \ti{n} (\lambda)+1.
    \]
Therefore it follows from Proposition~\ref{B-S}   that
  \[
 \ti{n}(\lambda)\leq N(\lambda)\leq \ti{n}(\lambda)+1.
 \]
 Since $\ti{n}(\lambda)\leq \| \wt{B} (\lambda)\|^2_2$, Lemma~\ref{BS1} implies that the total number $N=N(\pi)$ of eigenvalues of the operator $H$ lying above the point $\pi$ satisfies the bound
  \[
  N \leq \| \wt{B}(\pi)\|^2_2+1.
 \]
 
 According to  \e{eq:M4x} we have
 \[
 \| \wt{B}(\pi)\|_2\leq \pi^{-1} \big( \| V\|_2
 +  \| Q^\alpha V Q^\alpha \| \| Q^{-\alpha}  \wt{A} Q^{-\alpha}\|_2\big).
 \]
 Using \e{eq:M4y}, we find that
   \begin{equation}
 \| Q^{-\alpha}  \wt{A} Q^{-\alpha}\|_2^2= \pi^{-4}
 \int_{0}^\infty \int_{0}^\infty  \langle \ln t\rangle^{-2\alpha}(ts)^{-1} \sigma_{1} (t/s)^2 \langle \ln s\rangle^{-2\alpha} dtds=:\gamma_{\alpha}^2.
 \label{eq:L117}\end{equation}
 This yields the bound
   \begin{equation}
  N \leq \pi^{-2} \big( \| V\|_2
 +   \gamma_{\alpha} \| Q^\alpha V Q^\alpha \| \big)^2+1.
 \label{eq:L114}\end{equation}
 
 Observe further that according to representation \e{eq:M3} in the case $w\neq 0$, the operator $B(\lambda)$ has an eigenvalue which tends to $+\infty$ as $\lambda\to\pi$. Therefore Proposition~\ref{B-S} shows that the operator $H$ has at least one eigenvalue above the point $\lambda=\pi$. Note also that 
 \begin{equation}
    \| w\|^2= (V\psi_{0}, \psi_{0})= \int_{0}^\infty \int_{0}^\infty  {\bf v}(t,s)(ts)^{-1/2}  dtds.
 \label{eq:L115}\end{equation}
 In particular, in the Hankel case ${\bf v} (t,s)= v(t+s)$, we have
  \[
    \| w\|^2= \pi \int_{0}^\infty v (t)  dt. 
 \]
 Since a nonnegative Hankel operator $V$  has necessarily   nonnegative kernel $ v(t)$, the equality $w=0$ implies that $V=0$, and hence $w\neq 0$ if $V$ is nontrivial.

 Let us summarize the results obtained.

 \begin{theorem}\label{DS}
 Suppose that $V\geq 0$, $V\in {\goth S}_{2}$ and that inclusion \e{eq:H4c} is true for some $\alpha>3/2$. 
  Then   the total number $N $ of eigenvalues of the operator $H$ larger than $\pi$
  is finite and satisfies the bound \e{eq:L114} where the constant $\gamma_{\alpha}$ is defined by equalities \e{eq:L117} and \e{eq:L1s}. 
  
  Suppose additionally that integral \e{eq:L115} is not zero $($or that $V\neq 0$ in the Hankel case$)$. 
  Then     the operator $H$ has at least one eigenvalue larger than $\pi$. 
 \end{theorem}
 
 If $V$ has a negative part, then of course it can only diminish the   number $N $. In this case $V$ in the right-hand side of  \e{eq:L114} should be replaced by $V_{+}=(V+|V|)/2$.


\end{document}